\documentclass[12pt,a4paper,twoside]{article}
\usepackage{amsthm,amsmath,amssymb}
\newtheorem{teor}{Theorem}[section]
\newtheorem{defin}[teor]{Definition}
\newtheorem{lemm}[teor]{Lemma}
\newtheorem{osse}[teor]{Remark}
\newtheorem{prop}[teor]{Proposition}
\newtheorem{defi}[teor]{Definition}
\newtheorem{coro}[teor]{Corollary}
\newtheorem{prob}[teor]{Problem}

\newcommand{\bele}{\begin{lemm}\begin{sl}}
\newcommand{\enle}{\end{sl}\end{lemm}}
\newcommand{\bedef}{\begin{defi}\begin{sl}}
\newcommand{\eddef}{\end{sl}\end{defi}}
\newcommand{\bete}{\begin{teor}\begin{sl}}
\newcommand{\ente}{\end{sl}\end{teor}}
\newcommand{\beos}{\begin{osse}\begin{rm}}
\newcommand{\eddos}{\end{rm}\end{osse}}
\newcommand{\bepr}{\begin{prop}\begin{sl}}
\newcommand{\empr}{\end{sl}\end{prop}}
\newcommand{\bepro}{\begin{prob}\begin{rm}}
\newcommand{\empro}{\end{rm}\end{prob}}
\newcommand{\bede}{\begin{defin}\begin{sl}}
\newcommand{\edde}{\end{sl}\end{defin}}
\newcommand{\beco}{\begin{coro}\begin{sl}}
\newcommand{\enco}{\end{sl}\end{coro}}

\newcommand{\thspace}{\hspace{3mm}}

\newcommand{\quext}{\quad\text}
\newcommand{\qquext}{\qquad\text}

\newcommand{\RR}{\mathbb{R}}

\newcommand{\beeq}[1]{\begin{equation}\label{#1}}
\newcommand{\eddeq}{\end{equation}}

\newcommand{\beeqa}[1]{\begin{eqnarray}\label{#1}}
\newcommand{\eddeqa}{\end{eqnarray}}

\newcommand{\beal}[1]{\begin{align}\label{#1}}
\newcommand{\eddal}{\end{align}}

\newcommand{\bespl}[1]{\begin{split}\label{#1}}
\newcommand{\edspl}{\end{split}}

\newcommand{\bega}[1]{\begin{gather}\label{#1}}
\newcommand{\edga}{\end{gather}}

\newcommand{\beeqax}{\begin{eqnarray*}}
\newcommand{\eddeqax}{\end{eqnarray*}}

\def\qed{\ifmmode 
  \else \leavevmode\unskip\penalty9999 \hbox{}\nobreak\hfill
  \fi
  \quad\hbox{\hskip.5em\vrule width.4em height.6em depth.05em\hskip.1em}}
\def\endproofsym{\qed}
\renewenvironment{proof}[1][Proof]{\trivlist\item[\hskip\labelsep{\hskip0pt
    {\normalfont\scshape#1.}\hskip .321429\parindent}]\ignorespaces}
{\endproofsym\endtrivlist}
\def\endnobox{\def\endproofsym{}\end{proof}\def\endproofsym{\qed}}

\newcommand{\no}{\nonumber}

\newcommand{\beeqao}{\begin{eqnarray}\no}
\newcommand{\bealo}{\begin{align}\no}
\newcommand{\besplo}{\begin{split}\no}
\newcommand{\begao}{\begin{gather}\no}

\newcommand{\nor}[2]{\|#1\|_{#2}}

\newcommand{\+}{\hspace{1pt}}
\newcommand{\perogni}{\forall\,}
\newcommand{\esiste}{\exists\,}

\newcommand{\io}{\int_\Omega}

\newcommand{\epsi}{\varepsilon}

\newcommand{\rhs}{right hand side}

\DeclareMathOperator{\deriv}{d}

\newcommand{\LDtH}{L^2(0,t;H)}

\let\TeXchi\chi
\def\chi{{\setbox0 \hbox{\mathsurround0pt
$\TeXchi$}\hbox{\raise\dp0 \copy0 }}}

\newcommand{\teta}{\vartheta}

\newcommand{\J}{{\cal J}}

\newcommand{\zzu}{_{0,1}}
\newcommand{\zzd}{_{0,2}}

\newcommand{\barO}{\overline{\Omega}}

\newcommand{\dit}{\deriv\!t}

\newcommand{\dix}{\deriv\!x}
\newcommand{\diy}{\deriv\!y}

\newcommand{\ddt}{\frac{\deriv\!{}}{\dit}}

\numberwithin{equation}{section}
\begin{document}

\title{Finite-dimensional global attractor \\ for a nonlocal phase-field system
\thanks{
 This work was partially supported by the Italian MIUR-PRIN Research Project
2008 {\it Transizioni di fase, isteresi e scale multiple}}
}
\author{Maurizio Grasselli\\
{\sl Dipartimento di Matematica ``F.~Brioschi"}\\
{\sl Politecnico di Milano}\\
{\sl Via E.~Bonardi, 9}\\
{\sl I-20133 Milano, Italy}\\
{\tt maurizio.grasselli@polimi.it}}
\maketitle

\begin{abstract}\noindent
 We analyze a phase-field system where the energy balance equation
 is linearly coupled with a nonlinear and nonlocal ODE for the order parameter $\chi$.
 The latter equation is characterized by a space convolution term which models particle interaction
 and a singular configuration potential that forces $\chi$ to take values in $(-1,1)$.
 We prove that the corresponding dynamical system has a bounded absorbing set in a
 suitable phase space. Then we establish the existence
 of a finite-dimensional global attractor.
\end{abstract}

 \vspace{2mm}

\noindent
{\bf AMS (MOS) subject clas\-si\-fi\-ca\-tion:}\thspace
35B41, 35Q99, 80A22.

\vspace{2mm}



\pagestyle{myheadings}
\newcommand\testopari{\sc M. Grasselli}
\newcommand\testodispari{\sc attractor for a nonlocal phase-field system}
\markboth{\testodispari}{\testopari}


\section{Introduction}
\label{intro}

Consider a two-phase material (like, for instance, a mixture of ice and water)
which occupies a
bounded domain in $\Omega \subset \RR^d$, $1\le d\le3$ and
denote by $\teta$ its relative temperature with respect to a given
(constant) critical one (e.g., the one at which the two phases coexist).
A well-known model which accounts for phase changes due to the variation of $\teta$ only
was proposed and analyzed in \cite{Cag} (see also \cite{BS} and references therein).
This is based on the Ginzburg-Landau theory of phase transitions and it
assumes as further variable an order parameter (or phase-field) $\chi$ which characterizes,
for instance, the most energetic phase of the material (say water,
in a water-ice system). Using a phenomenological argument, it is postulated that
the evolution of $\chi$ is ruled by a {\sl gradient flow} of the form
\beeq{pargf}
  \chi_t=-\delta_\TeXchi E(\teta,\chi),
\end{equation}
$\delta_\TeXchi$ being the variational derivative with respect to
$\chi$ of the {\sl free energy} $E$ defined by
\beeq{energydiff}
   E(\teta,\chi)=\io\Big(\frac\nu2|\nabla\chi(x)|^2
   +W(\chi(x))
   -\alpha\teta(x)\+\chi(x)\Big)\,\dix,
\end{equation}
where $\nu>0$ and $\alpha\in\RR$ are given constants, the latter one being
related to the latent heat. Here $W$  is the (density of) potential
energy associated with the phase configuration which can be defined
either on a finite interval (e.g. $(-1,1)$) or on the whole real line. In
the first case, $W$ is called {\sl singular} and the most typical form is
the following
\begin{equation}
\label{loga}
W(r)=(1+r)\ln(1+r)
  +(1-r)\ln(1-r)-\frac{\lambda}{2} r^2,
  \quad r\in(-1,1),
\end{equation}
where $\lambda\in \RR$. However, though singular potentials are more acceptable from
the modeling viewpoint, $W$ is very often approximated by a
polynomial function defined on $\RR$ like, e.g., $W(r)=(r^2 -1)^2$.
In this case the potential is called {\sl smooth}.

It is worth observing that the term on the \rhs\ of \eqref{pargf} may
be viewed as a sort of {\sl generalized force} driving the evolution of
$\chi$, that is, the phase transformation. Combining now the balance
equation for the internal energy density with \eqref{pargf} we find
the evolution system
\beal{pf1}
  & (\teta +\alpha\chi)_t-\Delta\teta=f,\\
 \label{pf2}
  & \chi_t -\nu\Delta \chi + W'(\chi)-\alpha\teta=0,
\end{align}
in $\Omega\times (0,+\infty)$, where $f$ represents a given
(volumic) heat source and some constants have been taken equal to
one for simplicity. Equations \eqref{pf1}-\eqref{pf2} are also known
as Caginalp phase-field system. An important feature of this rather simple
model is the fact that its (formal) asymptotic limits are well-known sharp
interface problems (see, e.g., \cite{Cag2,CagC,CagF}). We recall that
phase-field systems are diffuse interface models in the sense that the
sharp interface separating two different phases, which is usually very
difficult to handle theoretically and numerically,
is replaced by the level set of a suitable (smooth) order parameter.

From the mathematical viewpoint, system \eqref{pf1}-\eqref{pf2}
endowed with initial conditions and various types of boundary
conditions has been investigated by many authors  (see,
e.g., \cite{BZ,BrCH,BrH,CGGM,CM,CFP,DKS,EZ,GMS,GPS2,JR,Kal,Kap,Lau,SA,Sc,Zha}
and references therein),
even under more general assumptions (for instance, where $\alpha$ is a function
depending on $\chi$). Besides well-posedness results, an important issue
is the longtime behavior of solutions. This behavior is usually non
trivial since the set of stationary states is a {\sl continuum} if
the spatial dimension is greater than one. Thus the existence of a Lyapunov
functional does not guarantee that a given trajectory converges to a single equilibrium.
To prove that the {\L}ojasiewicz-Simon approach has been employed (see, in particular,
\cite{CGGM,CM,CFP,GPS2,Zha}). Regarding the global dynamics, system
\eqref{pf1}-\eqref{pf2}, endowed with appropriate boundary
conditions, can be interpreted as a dynamical system in a suitable
phase space. This system is dissipative, i.e., there exists a bounded
absorbing set in the phase space. In addition, it possesses a
finite-dimensional global attractor as well as an exponential attractors
(cf., for instance, \cite{BZ,BrCH,BrH,CGGM,GMS,Kal,Kap,Lau}, see also
\cite{GGM,MZ} for related results).

On the other hand, the free energy $E$ can be viewed as an
approximation of a nonlocal expression of the following type (see
\cite{CF} and references therein, cf. also \cite{Ga,GZ,GL1,GL2} for
phase separation models)
\begin{equation}
\label{energynl} E_{\tt
nl}(\teta,\chi)=\io\Big(\io k(x-y)\frac{|\chi(x)-\chi(y)|^2}4
  \,\diy +W(\chi(x))
    -\alpha\teta(x)\+\chi(x)\Big)\,\dix,
\end{equation}
where $k:\RR^d \to \RR$ is an interaction kernel satisfying
$k(x)=k(-x)$ and such that
\begin{equation}
\label{kappa}
\kappa(x):=\io k(x-y)\,\diy
\end{equation}
is bounded and nonnegative (see, for instance, \cite[Rem.~2.2]{GZ} for concrete
examples). Indeed one can recover \eqref{energydiff} by
taking
\begin{equation*}
\io \frac\nu2|\nabla\chi(x)|^2\,\dix,
\end{equation*}
as a first approximation of
\begin{equation*}
\io \io k(x-y) |\chi(x)-\chi(y)|^2
  \,\diy \dix,
\end{equation*}
where $\nu= 2 \int_\Omega k(y)(y_i)^2 dy$ is supposed to be
independent of coordinate $i$. In concrete examples $k$ is localized in a
neighborhood of $0$ so that $\nu$ is related to the thickness of the
interface between the components. It
is interesting to point out that this kind of approximation was already
introduced by Van der Waals in his celebrated paper \cite{Ro} and since
then it was widely adopted in the mathematical literature on phase
transitions. This might be due to the fact that most people are more
used to deal with PDE rather than nonlocal operators.

The evolution system which derives from the nonlocal free energy
\eqref{energynl} takes the form \beal{calorein}
  & \teta_t+\alpha\chi_t-\Delta\teta=f,\\
 \label{phasein}
  & \chi_t+ \kappa\chi + W^\prime(\chi)=\J[\chi]+\alpha\teta,
\end{align}
in $\Omega\times (0,+\infty)$, where
\beeq{spaceconv}
  \J[v](x)=\io k(x-y)v(y)\,\diy, \qquad x\in\Omega.
\end{equation}
The well-posedness for this system was firstly established in
\cite{BCW} through a semigroup approach in the case
$\Omega=\RR$ for a smooth $W$. In a similar setting, existence of
traveling waves for small $\alpha$ was proven in \cite{BC1} (see
also \cite{BBH,BC2,BCh,BFRW,Fi,Wa} for results on the nonlocal
Allen-Cahn equation, i.e., \eqref{phasein} with $\alpha=0$). Nucleation
simulation by using nonlocal interactions has been studied in \cite{RRGE}. Phase
interface dynamics associated with \eqref{calorein}-\eqref{phasein}
was analyzed in \cite{CF} and \cite{CCE} (cf. also references therein)
by formal asymptotics (cf. also \cite{CE} for higher-order approximations and
\cite{D} for a related model). 
The case of bounded multi-dimensional domains was examined in \cite{BHZ} and
\cite{FIP2} (cf. also \cite{KS1,KS2,KRS1,KRS2,SZ} for results on
more refined models and \cite{ABH} for the numerical analysis).
Well-posedness issues were discussed in \cite{BHZ} when $\teta$ is
subject to homogeneous Neumann boundary condition and $W$ is
smooth. In addition, some results concerned with the asymptotic
behavior were demonstrated (e.g., the existence of a bounded
absorbing set in one spatial dimension). In \cite{FIP2}, the potential
$W$ is smooth as well, but $\teta$ is subject to Dirichlet
homogeneous. There the main goal was to establish the convergence
of a given trajectory to a single stationary state an this fact was
proven by means of a generalized version of the {\L}ojasiewicz-Simon
inequality, provided that $W$ is real analytic. Note that, also in the
nonlocal case, the set of stationary states can be a {\sl continuum} in
dimension greater than one, even though $z\mapsto
W^\prime(z)+\kappa z$ is invertible (see \cite[Introduction]{FIP2},
cf. also \cite{BCh} for the nonlocal Allen-Cahn equation). A similar
issue is analyzed in the more recent contribution \cite{GPS3}, where
equation \eqref{phasein} has an additional relaxation term of the
form $\varepsilon \chi_{tt}$ and $W$ is a singular potential (though
not of type \eqref{loga}).

However, none of the above results is concerned with the existence of
a global attractor even though this is a rather natural feature of
phase-field dynamics (cf. above for the local models). A possible explanation
might be that it is not evident how to show some precompactness for the variable
$\chi$. In this note we show how the existence of a finite-dimensional
global attractor can be established. This is achieved by estimating the
differences of two trajectories which originated from a suitable
bounded absorbing set. Such a difference is controlled by a contractive
part depending on the initial data plus a precompact term. For the
latter it is essential (but {\sl natural}) to assume some compactness
on the linear nonlocal operator $\J$. For the sake of simplicity, we will
take full advantage of the setting developed in \cite{GPS3}. Thus we
will suppose that $\teta$ satisfies Dirichlet homogeneous boundary
condition, while $W(r)$ will be a singular potential defined on $(-1,1)$
which goes to $+\infty$ as $r$ approaches the pure states $\pm1$.
This will allow us to easily deduce well-posedness  since the model
studied in \cite{GPS3} essentially reduces to ours when we take
$\varepsilon=0$. Nonetheless, we think that our argument can also be
applied to different boundary conditions for $\teta$ (e.g.,
homogeneous Neumann) and, especially, to other types of potentials
like, for instance, \eqref{loga} (see Remarks \ref{loga1} and
\ref{logpot} below).

The paper is organized as follows. In the next section we discuss
well-posedness issues. Section \ref{dissipative} contains a basic
uniform in time estimate which is exploited in Section \ref{globattr}
to show that our problem defines a dissipative dynamical system on a
suitable phase space. Then we conclude with our main result: the
existence of a finite-dimensional global attractor.


\section{Well-posedness}

\label{mainres}

Let us introduce some notation first. Set
$H:=L^2(\Omega)$ and $V:=H^1_0(\Omega)$ and denote by $(\cdot,\cdot)$
the scalar product in $H$ and by $\|\cdot\|$ the induced norm.
In general, $\nor\cdot X$ will indicate the norm in a generic real
Banach space $X$. Let $V$ be endowed with the norm
$\|\cdot\|_V:=\|\nabla\cdot\|$.
Let us identify $H$ with its topological dual $H'$ so that we have
the continuous and compact inclusions $V\subset H\subset V'$.
Moreover, we set $A:= -\Delta:D(A)=H^2(\Omega)\cap V \to H$.

Regarding the potential $W$, we suppose
\begin{align}
\label{W1}
&W\in C^2((-1,1);\RR^+),
   \qquad \lim_{r\to1^-,-1^+}W(r)=+\infty,\\
 \label{W2}
 &\esiste\lambda\in\RR:~~W''(r)\ge-\lambda, \qquad
  \perogni r\in  (-1,1).
\end{align}
It is easy to see that \eqref{W1}--\eqref{W2} entail
\beeq{W2.5}
  W'(r)r\ge W(r)-\frac{\lambda r^2}2-W(0),
   \quad\perogni r\in (-1,1).
\end{equation}
Moreover, if $v\in L^\infty(\Omega)$ is any function such that
$W(v)\in L^\infty(\Omega)$, then \beeq{sepprop}
  \esiste \delta=\delta\big(W,\nor{v}{L^\infty(\Omega)}\big)\in (0,1)
   \quext{such that }\,-1+\delta\le v(x)\le 1-\delta
   \quext{for a.e.~}\,x\in\Omega.
\end{equation}
Consider now $\kappa$ defined by \eqref{kappa} and the nonlocal operator $\J$
(cf.~\eqref{spaceconv}).
We require (cf., for instance, \cite{GZ})
\bega{J1}\tag{J1}
\int_{\Omega}\int_{\Omega} \vert k(x-y)\vert dydx =k_0 <+\infty,\\
 \label{J2}\tag{J2}
\textrm{ess} \sup_{x\in\Omega} \io \vert k(x-y)\vert \diy = k_1 <+\infty.
\end{gather}
Observe that the above assumptions entail that  $\J\in
\mathcal{L}(H;H)$ is self-adjoint and compact. Also, $\J$  is compact
from $L^\infty(\Omega)$ to $C(\barO)$.
\noindent%
System \eqref{calorein}--\eqref{phasein} with $\teta$ subject to
Dirichlet homogeneous boundary condition and initial conditions can
now be written as follows \beal{calore}
  & \teta_t+\alpha\chi_t+A\teta=f, \qquext{a.e.~in }\,\Omega\times(0,T),\\
  \label{phase}
  & \chi_t+ \kappa\chi+W'(\chi)=\J[\chi]+\alpha\teta, \qquext{a.e.~in }\,\Omega\times(0,T),\\
  \label{iniz1in}
  &\teta_{|t=0}=\teta_0,\quad \chi_{|t=0}=\chi_0\qquext{a.e.~in }\,\Omega,
\end{align}
where \begin{equation} \label{regteta}
  \teta_0\in V,
\end{equation}
and
\begin{equation}
\label{regf}
  f\in H.
\end{equation}
Moreover, we suppose that
\begin{equation}
\label{regchi}
  \chi_0\in L^\infty(\Omega)\quext{s.t. }
   \esiste \delta_0\in (0,1):~~-1+\delta_0\le\chi_0
         \le 1-\delta_0 \quext{a.e.~in }\,\Omega,
\end{equation}
where the latter property
is equivalent to say that $W(\chi_0)\in L^\infty(\Omega)$.

Well-posedness can be proven arguing as in the proof of \cite[Thm.2.1]{GPS3} with $\epsi=0$ (cf. also \cite[Rem.2.7]{GPS3})
\bete\label{teoexist}
 Let\/ \eqref{W1}--\eqref{W2}, \eqref{J1}--\eqref{J2},
 \eqref{regteta}--\eqref{regchi} hold.
 Then, for any given $T>0$, there exists one and only one
 pair $(\teta,\chi)$ such that
 \bega{regoteta}
    \teta\in L^2(0,T; D(A)) \cap H^1(0,T;H)\cap C^0([0,T];V),\\
  \label{regochi}
    \chi,~W(\chi)\in L^\infty(0,T;L^\infty(\Omega)),\quad
    \chi_{t} \in L^2(0,T;L^\infty(\Omega)),
 \end{gather}
 which solves \eqref{calore}-\eqref{iniz1in}.
 Moreover, there exists $\delta=\delta(W,k_0,\kappa, \alpha, f,\delta_0)\in (0,1)$ such that
 \begin{equation*}
 -1+\delta\le\chi(t)
         \le 1-\delta \quext{a.e.~in }\,\Omega,
 \end{equation*}
 for almost any $t\in(0,T)$.
%
%
 Next, given two triplets $(\teta\zzu,\chi\zzu)$,
 $(\teta\zzd,\chi\zzd)$ of initial data satisfying
 conditions\/ \eqref{regteta}--\eqref{regchi} (the latter w.r.t.~possibly
 different constants $\delta_{0i}>0$, $i=1,2$) and denoting
 the corresponding solutions
 by~$(\teta_1,\chi_1)$, $(\teta_2,\chi_2)$, respectively, we
 have the continuous dependence estimate
 \bealo
  & \|(\teta_1-\teta_2)(t)\|
   +\|\nabla(\teta_1-\teta_2)\|_{\LDtH}
   +\|(\chi_1-\chi_2)(t)\|\\
 \label{dipcon}
  & \mbox{}~~~~~
   \le \Lambda_0\big(\|\teta\zzu-\teta\zzd\|
    +\|\chi\zzu-\chi\zzd\|\big),
    \quad\perogni t\in[0,T],
 \end{align}
 where the positive constant $\Lambda_0$ depends on $T$, $\Omega$,
 $W$,  $k_0$, $k_1$, $\alpha$, $f$, and on the initial data (in particular, on
 $\delta_{0i}$, $i=1,2$).
\ente

\beos
To establish existence instead of the fixed-point technique used in \cite{GPS3}
one can use a vanishing viscosity argument like in \cite{FIP2}.
\eddos

We also have a higher-order control for the temperature difference, namely,

\beco
\label{addcontdep}
Let the assumptions of Theorem \ref{teoexist} hold. In addition to estimate
\eqref{dipcon} we have
\begin{equation}
 \label{dipconbis}
 \|\nabla(\teta_1-\teta_2)(t)\|
   +\|A(\teta_1-\teta_2)\|_{\LDtH}
    \le \Lambda_1\big(\|\teta\zzu-\teta\zzd\|_V
    +\|\chi\zzu-\chi\zzd\|\big),
\end{equation}
for all $t\in [0,T]$ where $\Lambda_1$ is a positive constant similar to $\Lambda_0$.
\enco

\begin{proof} It suffices to write equation \eqref{calore} for $\teta_1-\teta_2$,
multiplying by $A(\teta_1-\teta_2)$ recalling \eqref{phase} and using \eqref{dipcon}.
\end{proof}

\beos \label{kernel}
Assumptions \eqref{J1}-\eqref{J2} are satisfied
by concrete examples of interaction kernels like the ones mentioned in \cite[Rem.~2.2]{GZ}.
\eddos

\beos \label{weakhyp} Assumptions \eqref{regteta} and
\eqref{regf} can be weakened as follows
\begin{equation}
\teta \in D(A^{\frac{\rho}{2}}), \quad f \in D(A^{(\rho-1)/2}),
\end{equation}
where $\rho\in (\frac34,1)$. In this case, we recall that $D(A^\varrho)
\subset H^{2\varrho}(\Omega) \hookrightarrow L^\infty(\Omega)$ since the spatial dimension
is three at most.
In this case estimate \eqref{dipconbis} becomes
\begin{align*}
  & \|A^{\frac{\rho}{2}}(\teta_1-\teta_2)(t)\|
   +\|A^{\frac{\rho+1}{2}}(\teta_1-\teta_2)\|_{\LDtH}
  \\
  & \mbox{}~~~~~
   \le \Lambda_1\big(\|A^{\frac{\rho}{2}}(\teta\zzu-\teta\zzd)\|
    +\|\chi\zzu-\chi\zzd\|\big),
    \quad\perogni t\in[0,T].
 \end{align*}

\eddos
\beos\label{loga1}
 The second assumption \eqref{W1} does not cover, for instance,
 the case of a potential which is bounded in $[-1,1]$ like
 the logarithmic potential \eqref{loga}.
This important class of potentials can still be handled if we deal with
the Allen-Cahn equation only and temperature is assumed to be given
(see Remark \ref{logpot} below). However, system
\eqref{calorein}-\eqref{phasein} with more general potentials
require further arguments and will be analyzed elsewhere (see
\cite{GS}). \eddos


\section{A dissipative estimate}

\label{dissipative}

Here we establish some uniform in time estimates which are essentially
contained in the proof of \cite[Thm.2.1]{GPS3} taking $\epsi=0$.
\bete\label{teodiss} Let the assumptions of Theorem \ref{teoexist}
hold. Then the unique solution $(\teta,\chi)$ to
\eqref{calore}-\eqref{iniz1in} satisfies the following estimate
\begin{equation}
\label{diss1}
 \|\teta(t)\|^2_V + \Vert W(\cdot,\chi(t))\Vert_{L^\infty(\Omega)}
 \leq C_0 (1+ c_{\delta_0} + \Vert \teta_0\Vert^2_V)e^{-\beta t} + C_1,
\end{equation}
for any $t\geq 0$. Here $\beta$, $C_0$ and $C_1$ are positive constants
which depend on $\Omega$,
$\lambda$, $W(0)$, $k_1$, $\alpha$ and $f$
at most, while $c_{\delta_0}>0$ depends on
$\delta_0$ and $W$. \ente
\noindent%
\begin{proof}
Let us multiply equation \eqref{calore} by $\teta(t) + \xi A\teta(t)$
for some given $\xi>0$ to be chosen in the sequel. Integrating over
$\Omega$, we get
\begin{equation*}
\ddt\left(\frac12\Vert\teta\Vert^2 + \frac\xi 2 \Vert\nabla\teta\Vert^2\right)
    + \Vert\nabla \teta\Vert^2 + \xi \Vert A\teta\Vert^2=(f -\alpha\chi_t, \teta +\xi A\teta).
\end{equation*}
Consider now equation \eqref{phase}. Multiplying it by $\chi_t(t)+
\eta \chi(t)$, where $\eta>0$ will be fixed in the sequel, and integrating
over $\Omega$ we find
\begin{align*}
&\ddt\left(\frac\eta 2\Vert\chi\Vert^2 + \frac{1}{2}(\kappa\chi,\chi)+(W(\chi),1)\right)
    + \Vert\chi_t\Vert^2 + \eta(\kappa\chi,\chi)+\eta(W^\prime(\chi),\chi)\\
&=( \J[\chi]+\alpha\teta, \chi_t+\eta\chi).
\end{align*}
Adding the two identities we obtain
\begin{align*}
&\ddt \mathcal{E}
+ \Vert\nabla \teta\Vert^2 + \xi \Vert A\teta\Vert^2
+ \Vert\chi_t\Vert^2 + \eta(\kappa\chi,\chi)+\eta(W^\prime(\chi),\chi)
    \\
\nonumber
&=(f ,\teta +\xi A\teta)-\xi\alpha(\chi_t, A\teta)
+ ( \J[\chi], \chi_t+\eta\chi) +
\eta\alpha (\teta,\chi),
\end{align*}
where, for all $t\geq0$,
\begin{equation*}
\mathcal{E}(t)=
\frac12\Vert\teta(t)\Vert^2 + \frac\xi 2 \Vert\nabla\teta(t)\Vert^2 +
\frac\eta 2\Vert\chi(t)\Vert^2 + \frac{1}{2}(\kappa\chi(t),\chi(t))
+ (W(\chi(t)),1).
\end{equation*}
Recalling \eqref{W2.5}, \eqref{J1} and the fact that $\vert\chi\vert
\leq 1$ almost everywhere in $\Omega\times(0,T)$, for any given
$T>0$, it is not difficult to choose $\xi$ and $\eta$ such that
\begin{equation}
\label{diss2}
\ddt\mathcal{E} + c_1 \mathcal{E}
+ c_2\left(\Vert A\teta\Vert^2 + \Vert \chi_t\Vert^2\right)
\leq c_3 + c_4\Vert f \Vert^2,
\end{equation}
where $c_i$, $i=1,\dots,4$, are positive constants. In particular,
$c_1$ and $c_2$ depend on $\alpha$, while $c_3$ depends on
$\Omega$, $\lambda$, $W(0)$ and
$k_1$. Here we have also used the
Young and Poincar\'{e} inequalities.

We now argue as in \cite[Proof of Thm.2.1]{GPS3} and we  test
\eqref{phase} by $\chi_t+\sigma\chi$, for some $\sigma>0$ to be
properly selected, but we do not integrate over $\Omega$. We have
\begin{equation}
\ddt \mathcal{G}
    + (\chi_t)^2 + \sigma\kappa\chi^2 + \sigma W^\prime(\chi)\chi
    =( \J[\chi]+\alpha\teta)(\chi_t+\sigma\chi),
\end{equation}
where
\begin{equation*}
\mathcal{G}(x,t)=\frac\sigma 2(\chi(x,t))^2
+ \frac{1}{2}\kappa(x)(\chi(x,t))^2 + W(\chi(x,t)),
\qquad\text{ a.e. in }\Omega, \; t\geq 0.
\end{equation*}
Reasoning as before and using Young's inequality, we get
\begin{equation*}
\ddt \mathcal{G}(x,\cdot)
    + c_5 \mathcal{G}(x,\cdot) + c_6(\chi_t(x,\cdot))^2 \leq
    c_7(1+\vert \teta(x,\cdot)\vert^2), \qquad\text{ for a.a. }x\in\Omega,
\end{equation*}
where $c_i$, $i=5,\dots,7$ are positive constants depending at most
on $\Omega$, $\lambda$, $W(0)$, $k_1$ and $\alpha$.
Thanks to the continuous embedding
$D(A)\subset L^\infty(\Omega)$, we deduce
\begin{equation}
\ddt \mathcal{G}(x,\cdot)
    + c_5 \mathcal{G}(x,\cdot) + c_6(\chi_t(x,\cdot))^2 \leq
    c_8(1+\Vert A\teta\Vert^2),  \qquad\text{ for a.a. }x\in\Omega.
\end{equation}
If we multiply the above equation by $c_9=\frac{c_2}{2c_8}$ and
we add it to \eqref{diss2}, we obtain
\begin{equation}
\label{diss3}
\ddt(\mathcal{E} + c_9\mathcal{G}(x,\cdot))+ c_1 \mathcal{E} + c_9\mathcal{G}(x,\cdot)
+ \frac{c_2}{2}\Vert A\teta\Vert^2 + c_{10}\Vert \chi_t\Vert^2
\leq c_{11}(1+ \Vert f \Vert^2),
\end{equation}
for almost any $x\in\Omega$.

Applying now Gronwall's inequality to \eqref{diss3} we get, for all $t\geq 0$,
\begin{equation*}
\mathcal{E}(t)+ c_9\mathcal{G}(x,t)
\leq \left(\mathcal{E}(0) + c_9\mathcal{G}(x,0)\right) e^{-\mu t} +
\frac{2c_{11}}{\mu_1}\left(1+ \Vert f \Vert^2\right), \quad \text{ for a.a. }x\in\Omega,
\end{equation*}
where $\mu_1=\min\{1,c_1\}$, which yields \eqref{diss1}.
\end{proof}


\section{Existence of the global attractor}

\label{globattr}

A consequence of inequality \eqref{diss1} is that
$\mathcal{B}(R)=\{ u\in V \,:\, \Vert u\Vert \leq R\}$ for a fixed
$R>\sqrt{C_1}$ is absorbing for  $\teta(t)$ as well as for
$W(\chi(t))$. Therefore, consider, for instance,
\begin{equation}
\label{phasesp}
X= \mathcal{B}(R) \times \{v\in L^\infty(\Omega)\,:\,
\vert v\vert \leq 1-\delta_1, \;\text{ a.e. in }\Omega\},
\end{equation}
where $\delta_1 \in (0,1)$ is such that
$$
\{r\in (-1,1)\,:\, W(r) \leq
R^2\} \subseteq [-1+\delta_1,1+\delta_1].
$$
If, for each $(\teta_0,\chi_0)\in X$ and any $t\geq 0$, we define
$$
S(t)(\teta_0,\chi_0) = (\teta(t),\chi(t)),
$$
then, thanks to \eqref{diss1}, there exists $t_0>0$ such that $S(t)X\subseteq X$ for all $t\geq t_0$.
If we endow $X$ with the $V\times H$-metric then we obtain a
complete (bounded) metric space and $S(t)$ is strongly (Lipschitz)
continuous semigroup on $X$ owing to \eqref{dipcon} and
\eqref{dipconbis}. We now prove the main result of this note, namely,
the dynamical system $(X,S(t))$ has a finite-dimensional global
attractor.
\bete\label{teomain} Let the assumptions of Theorem \ref{teoexist} hold.
In addition, suppose that
\begin{equation}
\label{Wconv}
\lambda_0:= \textrm{ess}\inf_{x\in \Omega}\kappa(x) -\lambda >0.\\
\end{equation}
Then $(X,S(t))$ possesses a finite-dimensional connected global attractor.
\ente

\begin{proof} Consider $(\teta_{0i},\chi_{0i})\in X$, $i=1,2$,
set
\begin{equation*}
(\teta(t),\chi(t)) = ((\teta_1-\teta_2)(t),(\chi_1-\chi_2)(t))
\end{equation*}
where $(\teta_i(t),\chi_i(t))=S(t)(\teta_{0i},\chi_{0i})$ for $t\geq t_0$, and
observe that
\beal{calorediff}
  & \teta_t+\alpha\chi_t+A\teta=0, \qquext{a.e.~in }\,\Omega\times(t_0,+\infty),\\
  \label{phasediff}
  & \chi_t+\kappa\chi+W'(\chi_1(t)) - W'(\chi_2(t)) =\J[\chi]+\alpha\teta,
  \qquext{a.e.~in }\,\Omega\times(t_0,+\infty).
\end{align}
Let us multiply equation \eqref{calorediff} by $ A\teta(t)$.
Integrating over $\Omega$, we get
\begin{equation*}
\frac1 2 \ddt\Vert\nabla\teta\Vert^2
    +\Vert A\teta\Vert^2=-(\alpha\chi_t, A\teta).
\end{equation*}
from which, using the Young and Poincar\'{e} inequalities, we derive
the estimate
\begin{equation*}
\ddt\Vert\nabla\teta\Vert^2
    +c\Vert \nabla\teta\Vert^2 \leq c_\alpha \Vert\chi_t\Vert^2.
\end{equation*}
and, by comparison in \eqref{phasediff}, we deduce
\begin{equation}
\label{base0}
\ddt\Vert\nabla\teta\Vert^2
    +c\Vert \nabla\teta\Vert^2 \leq c(k_1,\delta_1,\alpha)\left (
    \Vert\chi\Vert^2 + \Vert\J[\chi]\Vert^2
    +\Vert\teta\Vert^2\right).
\end{equation}
for all $t\geq t_0$.
On the other hand, multiplying \eqref{phasediff} by $ \chi(t)$ and
integrating over $\Omega$ we find
\begin{equation*}
\frac12\ddt \Vert\chi\Vert^2 + (\kappa\chi,\chi)+ (W'(\chi_1) - W'(\chi_2),\chi)
    =( \J[\chi],\chi) +\alpha(\teta, \chi),
\end{equation*}
and \eqref{Wconv} entails
\begin{equation*}
\frac12\ddt \Vert\chi\Vert^2 + \lambda_0\Vert \chi\Vert^2
\leq ( \J[\chi],\chi) +\alpha(\teta, \chi).
\end{equation*}
Then, Young's inequality gives
\begin{equation*}
\frac12\ddt \Vert\chi\Vert^2 + \frac{\lambda_0}{2}\Vert \chi\Vert^2
\leq c_{\lambda_0} \Vert\J[\chi]\Vert^2 +c_\alpha\Vert\teta\Vert^2.
\end{equation*}
Since $\J$ is compact and self-adjoint we can find a finite-dimensional
projector $\Pi_{\lambda_0}$ such that
\begin{equation}
\label{proj}
\Vert\J[v]\Vert^2 \leq \frac{\lambda_0}{4c_{\lambda_0}}\Vert v\Vert^2 + \Vert \Pi_{\lambda_0}[v]\Vert^2,
\end{equation}
for all $v\in H$. As a consequence we have
\begin{equation}
\label{base1}
\frac12\ddt \Vert\chi\Vert^2 + \frac{\lambda_0}{4}\Vert \chi\Vert^2
\leq c_{\lambda_0} \Vert\Pi[\chi]\Vert^2 +c_\alpha\Vert\teta\Vert^2.
\end{equation}

Adding inequality\eqref{base0} multiplied by
$\mu_2=\frac{\lambda_0}{8c(k_1,\delta_1,\alpha)}$ to
\eqref{base1} and using \eqref{proj} yield
\begin{equation}
\label{base2}
\ddt \left(\mu_2\Vert\nabla\teta\Vert^2 + \Vert\chi\Vert^2 \right)
+c\mu_2\Vert\nabla\teta\Vert^2 + \frac{\lambda_0}{8}
\Vert \chi\Vert^2
\leq c(\lambda_0,k_1,\delta_1,\alpha)\left(\Vert \Pi[\chi]\Vert^2 +
\Vert\teta\Vert^2\right).
\end{equation}
Therefore, from \eqref{base2}, we deduce
\begin{align}
\label{base3}
\Vert\teta(t)\Vert_V^2 +\Vert \chi(t)\Vert^2 &\leq
c(\lambda_0,k_1,\delta_1,\alpha)e^{-\mu_3(t-t_0)}
\left(\Vert\teta(t_0)\Vert_V^2 + \Vert\chi(t_0)\Vert^2\right)\\
\nonumber
&+ c(\lambda_0,k_1,\delta_1,\alpha)\int_{t_0}^t \,\left(\Vert\teta(\tau)\Vert^2
+ \Vert \Pi_{\lambda_0}[\chi(\tau)]\Vert^2 \right) d\tau,
\end{align}
for all $t\in [t_0,T]$ and any fixed $T>t_0$. Here $\mu_3$ is a
positive constant depending on $\lambda_0$, $k_1$, $\delta_1$,
$\alpha$.

We now introduce the following pseudometric in $X$
$$
\mathbf{d}_T ((\teta_{01},\chi_{01}),(\teta_{02},\chi_{02}) )
= \left(\int_{t_0}^T\,\left(\Vert
(\teta_1 -\teta_2)(\tau)\Vert^2 + \Vert \Pi_{\lambda_0}[(\chi_1-\chi_2)(\tau)]
\Vert^2\right) d\tau\right)^{1/2}
$$
and we recall that a pseudometric is (pre)compact in $X$ (with
respect to the topology induced by the  $X$-metric) if any bounded
sequence in $X$ contains a Cauchy subsequence with respect to
$\mathbf{d}_T$ (see, for instance, \cite[Def.~2.57]{KM}).

Let us prove that $\mathbf{d}_T$ is precompact in $X$. Let
$\{(\teta_{0n},\chi_{0n})\}\subset X$ ($X$ is bounded) and set
$(\teta_n(t),\chi_n(t))=S(t)(\teta_{0n},\chi_{0n})$. Thanks to
\eqref{regoteta}, we have that $\{\teta_n\}$ is bounded in
$L^2(t_0,T;D(A))\cap H^1(t_0,T;H)$. Therefore it contains a
subsequence which strongly converges in $L^2(t_0,T,V)$. On the
other hand, we have that $\{\Pi_{\lambda_0}[\chi_n]\}$ is
bounded in $L^\infty(t_0,T;H)$. Also, by comparison in
\eqref{phase}, we deduce that $\{(\chi_n)_t\}$ is bounded in
$L^\infty(t_0,T;H)$. Therefore $\{(\Pi_{\lambda_0}[\chi_n])_t\}$
is bounded in $L^\infty(t_0,T;H)$ as well. Then
$\{\Pi_{\lambda_0}[\chi_{n}(\cdot)]\}$ contains a subsequence
which strongly converges in $L^2(t_0,T;H)$. Summing up
$\{(\teta_{0n},\chi_{0n})\}$ contains a Cauchy subsequence with
respect to $\mathbf{d}_T$.

From \eqref{base3}, we deduce that there exists $t^*>t_0$ such that
\begin{align*}
&\Vert S(t^*)(\teta_{01},\chi_{01}) - S(t^*)(\teta_{02},\chi_{02})
\Vert_X\\
&\leq
\frac{1}{2}\Vert(\teta_{01} - \teta_{02},\chi_{01} - \chi_{02})\Vert_X
+ C(\lambda_0,k_1,\delta_1,\alpha)
\mathbf{d}_{t^*} ((\teta_{01},\chi_{01}),
(\teta_{02},\chi_{02})).
\end{align*}
Hence $S(t)$ has a (connected) global attractor (see
\cite[Thm.~2.56, Prop.~2.59]{KM}) of finite fractal (i.e., box
counting) dimension (cfr. \cite[Thm.~2.8.1]{Ha}).
\end{proof}

\beos \label{convexity} Assumption \eqref{Wconv} seems
unavoidable when one wants to investigate the long-time behavior of
solutions (cf. \cite[(A4)]{BHZ} and \cite[(1.19)]{FIP2}).
In particular, thanks to this assumption,
if we take $f\equiv0$ and suppose $W$ real analytic then, on account of \cite[Thm.~2.6]{GPS3},
we have that the $\omega$-limit set of any pair $(\theta_0,\chi_0)$ satisfying
\eqref{regteta} and \eqref{regchi} reduces to a singleton $\{(0,\chi_\infty)\}$,
where
$$
\kappa\chi_\infty+W'(\chi_\infty)=\J[\chi_\infty], \qquad\text{a.e. in}\,\Omega.
$$
\eddos

\beos \label{weakphase} On account of Remark \ref{weakhyp},  we
could take a larger phase space by replacing $V$ with
$V_\rho=D(A^{\frac\rho2})$, $\rho\in (\frac34,1)$, in the
definition of $\mathcal{B}(R)$ and endowing $X$ with
$V_\rho\times H$-norm.\eddos

\beos \label{logpot} Consider the following equation
\begin{equation*}
\chi_t+\kappa\chi+W'(\chi)=\J[\chi]+g,
\end{equation*}
where $g\in L^\infty(\Omega\times(0,+\infty))$. In this case it is
possible to show a (uniform) separation property even when $W$ is a
more general potential like, e.g., \eqref{loga}. However, one should
use a comparison argument like, e.g., in \cite{GPS2}. Indeed, it is no
longer sufficient to show the global uniform boundedness of
$W(\chi)$ (cf. \eqref{diss1}). Then one can define a phase space
given by the second component of $X$ (see \eqref{phasesp}) and,
arguing as above for the $\chi$ component only, prove the existence
of a finite-dimensional global attractor.

\eddos

\bigskip\noindent
{\bf Acknowledgments.} The author thanks S.~Frigeri, G.~Schimperna, and
D.~Pra\v z\'ak for their helpful remarks. This work was partially supported
by the Italian MIUR-PRIN Research Project 2008 ``Transizioni di fase,
isteresi e scale multiple''.



\begin{thebibliography}{99}



\bibitem{ABH}
 S.~Armstrong, S.~Brown, J.~Han, {\sl Numerical analysis for a nonlocal phase field system},
 Int. J. Numer. Anal. Model. Ser.~B {\bf 1} (2010), 1--19.

\bibitem{BBH}
 P.W.~Bates, S.~Brown, J.~Han, {\sl Numerical analysis for a nonlocal Allen-Cahn equation},
 Int. J. Numer. Anal. Model. {\bf 6} (2009), 33--49.

\bibitem{BC1}
 P.W.~Bates, F.~Chen,
 {\sl Traveling wave solutions for a nonlocal phase-field system},
 Interfaces Free Bound.,
 {\bf 4} (2002),
 227--238.

 \bibitem{BC2}
 P.W.~Bates, F.~Chen,
 {\sl Spectral analysis and multidimensional stability of traveling waves
 for nonlocal Allen-Cahn equation},
 J.~Math.\ Anal.\ Appl.,
 {\bf 273} (2002),
 45--57.

\bibitem{BCW} P.W.~Bates, F.~Chen, J.~Wang, {\sl Global existence and uniqueness of solutions to a nonlocal phase-field system}, in: P.W. Bates, S.-N. Chow, K. Lu and X. Pan, Editors, US-Chinese Conference on Differential Equations and Applications, International Press, Cambridge, MA (1997), 14--21.

\bibitem{BCh} P.W.~ Bates, A.~Chmaj, {\sl An integrodifferential model
    for phase transitions:  stationary solutions in higher space
    dimensions}, J.~Statist.\ Phys., {\bf 95} (1999), 1119--1139.

\bibitem{BFRW}
P.W.~Bates, P.C.~Fife, X.~Ren, X.~Wang,
{\sl Traveling waves in a convolution model for phase transitions},
Arch.\ Rational Mech.\ Anal., {\bf 138} (1997), 105--136.




\bibitem{BHZ}
P.W.~Bates, J.~Han, G.~Zhao, {\sl
On a nonlocal phase-field system}, Nonlinear Anal., {\bf 64} (2006), 2251--2278.

\bibitem{BZ} P.W.~Bates, S.~Zheng, {\sl Inertial manifolds and
    inertial sets for the phase-field equations}, J.\ Dynamics Differential
    Equations, {\bf 4} (1992), 375--397.
    
\bibitem{BrCH}
D.~Brochet, X.~Chen, D.~Hilhorst, {\sl Finite dimensional exponential attractor for the phase field
model}, Appl.\ Anal., {\bf 49} (1993), 197--212.

\bibitem{BrH}
D.~Brochet, D.~Hilhorst, {\sl Universal attractor and inertial sets for the phase-field model},
Appl.\ Math.\ Lett., {\bf 4} (1991), 59--62.

\bibitem{BS}
 M.~Brokate, J.~Sprekels,
 ``Hysteresis and Phase Transitions'',
 Springer,
 New York,
 1996.

\bibitem{Cag}
 G.~Caginalp,
 {\sl An analysis of a phase field model of a free boundary},
 Arch.\ Rational Mech.\ Anal.,
 {\bf 92} (1986),
 205--245.

\bibitem{Cag2}
G.~Caginalp, {\sl Stefan and Hele-Shaw type models as asymptotic limits of the phase-field equations}, Phys.\ Rev.\ A, {\bf 39} (1989), 5887--5896.

\bibitem{CagC}
G.~Caginalp, X.~Chen, {\sl Convergence to the phase field model to its sharp interface limits},
European J.\ Appl.\ Math., {\bf 9} (1998), 417--445.

\bibitem{CE} G.~Caginalp, E.~Esenturk, {\sl A phase field
    model with non-local and anisotropic potential},
    Discrete Contin.\ Dyn.\ Syst.\ Ser.\ S,
    {\bf 4} (2011) 311--350.

\bibitem{CagF}
G.~Caginalp, P.C.~Fife, {\sl Dynamics of layered interfaces arising from phase boundaries},
SIAM J. Appl. Math., {\bf 48} (1988), 506--518.

\bibitem{CGGM}
C.~Cavaterra, C.G.~Gal, M.~Grasselli, A.~Miranville, {\sl Phase-field systems with nonlinear coupling and dynamic boundary conditions}, Nonlinear Anal., {\bf 72} (2010), 2375--2399.


\bibitem{CF} C.K.~Chen, P.C.~Fife, {\sl Nonlocal models of phase
    transitions in solids},  Adv.\ Math.\ Sci.\ Appl.,
    {\bf 10} (2000),
    821--849.
    
\bibitem{CCE} X.~Chen, G.~Caginalp, E.~Esenturk, {\sl A phase field
    model with non-local and anisotropic potential},
    Modelling Simul.\ Mater.\ Sci.\ Eng.,
    {\bf 19} (2011) 045006(8).

\bibitem{CM}
L.~Cherfils, A.~Miranville, {\sl Some results on the asymptotic behavior of the Caginalp
system with singular potentials}, Adv.\ Math.\ Sci.\ Appl., {\bf 16} (2007), 107--129.

\bibitem{CFP}
R.~Chill, E. Fa\v{s}angov\'{a}, J. Pr\"{u}ss, {\sl Convergence to steady states of solutions of the Cahn-
Hilliard and Caginalp equations with dynamic boundary conditions}, Math.\ Nachr., {\bf 13} (2006),
1448--1462.


\bibitem{DKS}
A.~Damlamian, N.~Kenmochi, N.~Sato, {\sl Subdifferential operator approach to a class of
nonlinear systems for Stefan problems with phase relaxation}, Nonlinear Anal., {\bf 23} (1994),
115--142.

\bibitem{D} N.~Dirr,
{\sl A Stefan problem with surface tension as the sharp interface limit of a nonlocal system of phase-field type},
J.~Statist.\ Phys., {\bf 114} (2004), 1085--1113.


\bibitem{EZ}
C.M.~Elliott, S.~Zheng, {\sl Global existence and stability of solutions to the phase-field
equations}, in "Free boundary problems," Internat.\ Ser.\ Numer.\ Math. {\bf 95}, 46--58,
Birkh\"{a}user Verlag, Basel, 1990.


\bibitem{FIP2}
 E.~Feireisl, F.~Issard-Roch, H.~Petzeltov\'{a},
 {\sl A non-smooth version of the {\L}ojasie\-wicz-Simon theorem with applications to
 non-local phase-field systems},
 J.\ Differential Equations,
 {\bf 199} (2004),
 1--21.

\bibitem{Fi}
P.C.~Fife, {\sl Well-posedness issues for models of phase transitions with weak interaction},
Nonlinearity,
{\bf 14} (2001), 221–-238.

\bibitem{Ga} H.~Gajewski, {\sl On a nonlocal model of
    non-isothermal phase separation},
    Adv.\ Math.\ Sci.\ Appl., {\bf 12} (2002),
    569--586.

\bibitem{GZ}
H.~Gajewski, K.~Zacharias,
{\sl On a nonlocal phase separation model},
J.\ Math.\ Anal.\ Appl.,
{\bf 286}  (2003),
11--31.

\bibitem{GGM}
C.G.~Gal, M.~Grasselli, A.~Miranville, {\sl Robust exponential attractors for singularly perturbed
phase-field equations with dynamic boundary conditions}, NoDEA Nonlinear Differential Equations Appl.,
{\bf 15} (2008), 535–-556.


\bibitem{GL1}
 G.~Giacomin, J.L.~Lebowitz,
 {\sl Phase segregation dynamics in particle systems with long range
   interactions. I. Macroscopic limits},
 J.~Statist.\ Phys.,
 {\bf 87} (1997),
 37--61.

\bibitem{GL2}
 G.~Giacomin, J.L.~Lebowitz,
 {\sl Phase segregation dynamics in particle systems with long range
   interactions. II. Interface motion},
 SIAM J.~Appl.\ Math.,
 {\bf 58} (1998),
 1707--1729.


\bibitem{GMS}
M.~Grasselli, A.~Miranville, G.~Schimperna, {\sl
The Caginalp phase-field system with coupled dynamic boundary conditions and singular potentials},
Discrete Contin.\ Dyn.\ Syst., {\bf 28} (2010), 67--98.

\bibitem{GPS2}
 M.~Grasselli, H.~Petzeltov\'a, G.~Schimperna,
 {\sl Long time behavior of solutions to the Caginalp system with
 singular potential},
 Z.\ Anal.\ Anwendungen,
 {\bf 25} (2006),
 51--72.

\bibitem{GPS3}
M.~Grasselli, H.~Petzeltov\'a, G.~Schimperna,
{\sl  A nonlocal phase-field system with inertial term,}
Quart.\ Appl.\ Math.,
{\bf 65} (2007), 451--469.

\bibitem{GS} M.~Grasselli, G.~Schimperna, in preparation.

\bibitem{Ha} J.K.~Hale, ``Asymptotic behaviour of
    dissipative systems", Amer. Math. Soc., Providence,
    RI, 1988.

\bibitem{JR}
A.~Jim\'{e}nez-Casas, A. Rodr\'{\i}guez-Bernal, {\sl Asymptotic behaviour for a phase field model
in higher order Sobolev spaces}, Rev.\ Mat.\ Complut., {\bf 15} (2002), 213--248.

\bibitem{Kal}
V.K.~Kalantarov, {\sl On the minimal global attractor of a system of phase field equations}, (Russian),
Zap. Nauchn. Sem. Leningrad. Otdel. Mat. Inst. Steklov. (LOMI) {\bf 188} (1991), Kraev.\ Zadachi
Mat.\ Fiz.\ i Smezh.\ Voprosy Teor.\ Funktsii. {\bf 22}, 70--86, 186 [translation in J.\ Math.\ Sci., {\bf 70}
(1994), 1767--1777].

\bibitem{Kap}
O.V.~Kapustyan, {\sl An attractor of a semiflow generated by a system of phase-field equations
without uniqueness of the solution} (Ukrainian), Ukra\"{\i}n.\ Mat.\ Zh., {\bf 51} (1999), 1006--1009 [Translation
in Ukrainian Math.\ J., {\bf 51} (1999), 1135--1139 (2000)].

\bibitem{KM}
N.J.~Koksch, A.J.~Milani, ``An introduction to semiflows", Chapman \& Hall/CRC, Boca Raton, FL, 2005.

\bibitem{KS1}
 P.~Krej\v c\'\i, J.~Sprekels,
 {\sl Nonlocal phase-field models for non-isothermal phase
   transitions and hysteresis},
 Adv.\ Math.\ Sci.\ Appl.,
 {\bf 14} (2004),
 593--612.

\bibitem{KS2}
 P.~Krej\v c\'\i, J.~Sprekels,
 {\sl Long time behavior of a singular phase transition model},
 Discrete Contin.\ Dyn.\ Syst.,
 {\bf 15} (2006), 1119--1135.

\bibitem{KRS1}
 P.~Krej\v c\'\i, E.~Rocca, J.~Sprekels,
 {\sl Nonlocal temperature-dependent phase-field models
   for non-isothermal phase transitions},
J.~Lond.\ Math.\ Soc. (2), {\bf 76} (2007), 197--210.

\bibitem{KRS2}
P.~Krej\v c\'\i, E.~Rocca, J.~Sprekels, {\sl
A nonlocal phase-field model with nonconstant specific heat},
Interfaces Free Bound., {\bf 9} (2007), 285--306.

\bibitem{Lau}
Ph.~Lauren\c{c}ot, {\sl Long-time behaviour for a model of phase-field type},
Proc.\ Roy.\ Soc.\ Edinburgh Sect.~A, {\bf 126} (1996), 167--185.

\bibitem{MZ}
A.~Miranville, S.~Zelik, {\sl Robust exponential attractors for singularly perturbed phase-field type
equations}, Electron.\ J.\ Differential Equations, {\bf 63} (2002), 1--28.

\bibitem{Ro}
J.S.~Rowlinson, Translation of J.D.~van der Waals, {\sl
The thermodynamic theory of capillarity under the hypothesis of a continuous variation of density},
J.\ Statist.\ Phys., {\bf 20} (1979), 197--244.

\bibitem{RRGE}
A.~Roy, J.M.~Rickman, J.D.~Gunton, K.R.~Elder, {\sl Simulation study of nucleation in 
a phase-field model with nonlocal interactions}, Phys.\ Rev.\ E, {\bf 57} (1998), 2610(8).

\bibitem{SA}
N.~Sato, T.~Aiki, {\sl Phase field equations with constraints under nonlinear dynamic boundary
conditions}, Commun.\ Appl.\ Anal., {\bf 5} (2001), 215--234.

\bibitem{Sc}
G.~Schimperna, {\sl Abstract approach to evolution equations of phase field type and applications},
J. Differential Equations, {\bf 164} (2000), 395--430.

\bibitem{SZ}
 J.~Sprekels, S.~Zheng,
 {\sl Global existence and asymptotic behaviour for a nonlocal phase-field model for
 non-isothermal phase transitions},
 J.\ Math.\ Anal.\ Appl.,
 {\bf 279} (2003),
 97--110.

\bibitem{Wa}
 X.~Wang,
 {\sl Metastability and stability of patterns in a convolution model
 for phase transitions},
 J.~Differential Equations,
 {\bf 183} (2002), 
 434--461.

\bibitem{Zha}
Z.~Zhang, {\sl Asymptotic behavior of solutions to the phase-field equations with Neumann boundary
conditions}, Commun.\ Pure Appl.\ Anal., {\bf 4} (2005), 683--693.



\end{thebibliography}
\end{document}